\documentclass[a4paper,11pt]{article}
\usepackage{}
\usepackage{mathrsfs}
\usepackage{mathrsfs}
\usepackage{mathrsfs}
\usepackage{bbm}
\usepackage{amssymb}
\usepackage{amsmath}
\usepackage{amsfonts}
\usepackage{amsmath}
\usepackage{graphicx}
\numberwithin{equation}{section}
\newtheorem{theorem2}{Theorem}[section]
\newtheorem{lemma2}{Lemma}[section]

\newtheorem {theorem3}{Theorem}[section]
\newtheorem {remark3}{Remark}[section]
\newtheorem {lemma3}{Lemma}[section]
\newtheorem {corollary3}{Corollary}[section]

\newtheorem {theorem4}{Theorem}[section]

\newtheorem {lemma4}{Lemma}[section]
\newtheorem {definition4}{Definition}[section]

\begin{document}
\title{\textbf{The blow-up phenomena and  exponential
decay of solutions  for a three-component
Camassa--Holm equations}}
\author{Xinglong Wu \footnote{E-mail: wxl8758669@aliyun.com}\\
Wuhan Institute of Physics and Mathematics,\\ Chinese Academy of Sciences,
 Wuhan 430071, P. R. China \smallskip\\
}
\date{}
\maketitle

\begin{abstract}
The present paper is mainly concerned with the blow-up phenomena and exponential decay of solution for a three-component
Camassa--Holm equation. Comparing with the result of Hu, ect. in the paper \cite{HLJ}, a new wave-breaking solution is obtained. The results of exponential decay of solution in our paper cover and extent the corresponding results in \cite{C-E1, H-M, W-G} .

\end{abstract}

\section{Introduction}
In this paper, we devote to the study of the Cauchy problem for a three-component
Camassa--Holm equation
\begin{equation}
\left\{\begin{array}{ll}m_t-m_{x}u+2mu_x+(mv+mw)_x+nv_x+lw_x=0,\\

n_t-n_{x}v+2nv_x+(nu+nw)_x+mu_x+lw_x=0,\\
l_t-l_{x}w+2lw_x+(lu+lv)_x+mu_x+nv_x=0,
\end{array}\right.
\end{equation}
which was introduced  by Qu and Fu in \cite{Q-F} to study multipeakons, where the potential $m=u-u_{xx}$, $n=v-v_{xx}$ and $l=w-w_{xx}$, $(t,x)\in\mathbb{R}^{+}\times\mathbb{R}$, and  the subscripts denote the
partial derivatives. The two peakon solitions of system (1.1) have the following form
\begin{equation*}
u(x,t)=p_1(t)\exp(-|x-q_1(t)|)+p_2(t)\exp(-|x-q_2(t)|),
\end{equation*}
\begin{equation*}
v(x,t)=r_1(t)\exp(-|x-q_1(t)|)+r_2(t)\exp(-|x-q_2(t)|),
\end{equation*}
\begin{equation*}
w(x,t)=s_1(t)\exp(-|x-q_1(t)|)+s_2(t)\exp(-|x-q_2(t)|),
\end{equation*}
where $p_i$,\;$q_i(t),\;r_i(t)$ and $s_i(t),i=1,2$ are functions of $t$, and the corresponding dynamical system  defined in \cite{Q-F}.
\par
Let the potential $v=w=0$, system (1.1) becomes
 the  classical Camassa-Holm (CH) equation in form
\begin{equation}
m_{t}+um_{x}+2u_{x}m=0, \qquad m=u-u_{xx}.
\end{equation}
which comes from an asymptotic approximation to the
Hamiltonian for the Green--Naghdi equations in shallow water theory. The CH equation models the unidirectional propagation of
shallow water waves over a flat bottom \cite{C-H}, and also is a
model for the propagation of axially symmetric waves in hyperelastic
rods \cite{Da}. It has a bi-Hamiltonian structure \cite{F-F} and is completely integrable
\cite{Co}, and with a Lax pair based on a linear spectral problem of
second order. Also, there exists smooth soliton solutions
on a non-zero constant background \cite{C-H-H}.
Compared with KdV equation, the CH equation not only approximates
unidirectional fluid flow in Euler's equations \cite{G-N} at the
next order beyond the KdV equation \cite{Ka1, Ka2}, but also there exists
blow-up phenomena of the strong solution and global existence of
strong solution \cite{Co, C-E, C-E1, C-E2}. It is remarkable that the
CH equation has peaked solitons of the form
$u(t,x)=ce^{-|x-ct|},\;c\in\mathbb{R}$ \cite{C-H-H}, which are
orbital stable \cite{C-S}, and $n$-peakon solutions \cite{B-S-S}
$$u(t,x)=\sum_{j=1}^{n}p_{j}(t)\exp(-|x-q_{j}(t)|),$$
where the positions $q_{j}$ and amplitudes $p_{j}$ satisfy the
system of ODEs
$$\left\{\begin{array}{ll}\dot{q_{j}}=\sum_{k=1}^{n}
p_{k}\exp(-|q_{j}-q_{k}|),\\
\\ \dot{p_{j}}=
p_{j}\sum_{k=1}^{n}p_{k}\text{sgn}(q_{j}-q_{k})\exp(-|q_{j}-q_{k}|),
\end{array}\right.$$
where $j=1,\cdots,n.$ The CH equation has attracted a lot
of interest in the past twenty years for various reasons \cite{B-C,
B-C1, Co1, C-E, C-S, Fu, W-Y}.
\par
If we neglect $w$ in system (1.1) to obtain 2-component CH equation, which is studied in \cite{F-Q}. They establish the local well-posedness in $H^s\times H^s,s>\frac{3}{2}$. Also, it has blow-up phenomena, if the initial data satisfy some condition.

Recently, Hu, Lin and Jin investigate the Cauchy problem of the three Camassa-Holm equation
(1.1) in \cite{HLJ} on the line. The authors
establish the local well-posedness, derive precise blow-up
scenario and the conservation law. Moreover, by the conservation law,  if the derivative of initial data is negative, they obtain the existence of strong solutions which blows up in finite time and derive the blow-up rate. In this paper, we give a new blow-up phenomena to system (1.1), as the initial data satisfy $$\int_{\mathbb{R}}(u_{0x}+v_{0x}+w_{0x})^3dx<-9E_0\sqrt{2E_0}.$$ Next, we study the exponential decay of the solution provided the initial data $z_0(x)=(u_0,v_0,w_0)\sim\mathcal {O}(e^{-\alpha|x|}),\;\alpha\in(0,1)$  as $x\rightarrow\pm\infty$ or the initial potential $(m_0,n_0,l_0)\sim \mathcal {O}(e^{-(1+\lambda)|x|}),\;\lambda>0$ as $x\rightarrow\pm\infty$. Moreover,
we get a class of traveling wave solutions of system (1.1).
\par
The remainder of the paper is organized as follows. In Section 2, we give a new wave-breaking
solution of system (1.1). In Section 3, the exponential decay of solution is established, if the initial data satisfy some decay condition. In Section 4,
we prove that system (1.1) has a class of traveling wave solution.\\
 Notation: For
simplicity, we identify all spaces of functions with function spaces
over $\mathbb{R}$, we drop $\mathbb{R}$ from our notation. For
$1\leq p\leq \infty$, the norm in the Banach space $L^p(\mathbb{R})$
will be written by $\|\cdot\|_{L^p}$, while $\|\cdot\|_{H^s},
s\in\mathbb{R}$ will stand for the norm in the classical Sobolev
spaces $H^s(\mathbb{R})$. We shall say for some $K>0$ that
$$f(x)\sim\mathcal {O}(e^{\alpha x})\quad\text{as}\;x\uparrow \infty,\quad \text{if}\quad
\lim_{x\rightarrow\infty}\frac{|f(x)|}{e^{\alpha x}}\leq K,$$ and
$$f(x)\sim o(e^{\alpha x})\quad\text{as}\;x\uparrow \infty,\quad \text{if}\quad
\lim_{x\rightarrow\infty}\frac{|f(x)|}{e^{\alpha x}}=0.$$
\section{The blow-up phenomena of solution}
\par
With the potential  $m=u-u_{xx}$, $n=v-v_{xx}$ and $l=w-w_{xx}$.
It is convenient to rewrite the system (1.1) in its formally equivalent differential
form
\begin{equation}
\left\{\begin{array}{ll}
 u_t+(u+v+w)u_x+(1-\partial_x^2)^{-1}(uv_x+uw_x)+\partial_x(1-\partial_x^2)^{-1}f=0,\\
v_t+(u+v+w)v_x+(1-\partial_x^2)^{-1}(vu_x+vw_x)+\partial_x(1-\partial_x^2)^{-1}g=0,\\
w_t+(u+v+w)w_x+(1-\partial_x^2)^{-1}(wu_x+wv_x)+\partial_x(1-\partial_x^2)^{-1}h=0,\\
 (u,v,w)|_{t=0}=(u_{0}(x),v_0(x),w_0(x)),\end{array}\right.
\end{equation}
where the functions $f,g$ and $h$ satisfy
\begin{equation}\begin{split}
&f(x,t)=u^2+\frac{1}{2}u_x^2+u_xv_x+u_xw_x+\frac{1}{2}v^2-\frac{1}{2}v_x^2+\frac{1}{2}w^2-\frac{1}{2}w_x^2,\\
&g(x,t)=v^2+\frac{1}{2}v_x^2+u_xv_x+w_xv_x+\frac{1}{2}u^2-\frac{1}{2}u_x^2+\frac{1}{2}w^2-\frac{1}{2}w_x^2,\\
&h(x,t)=w^2+\frac{1}{2}w_x^2+u_xw_x+v_xw_x+\frac{1}{2}u^2-\frac{1}{2}u_x^2+\frac{1}{2}v^2-\frac{1}{2}v_x^2.
\end{split}\end{equation}
Note that if choosing the Green function $G(x)=\frac{1}{2}e^{-|x|}, x\in \mathbb{R}, $ we have
$(1-\partial^{2}_{x})^{-1}f=G\ast f $ for all $ f\in
L^{2}(\mathbb{R})$. Then  Eq.(2.1) can  been rewritten  as follows
\begin{equation}
\left\{\begin{array}{ll}
 u_t+(u+v+w)u_x+G\ast(uv_x+uw_x)+\partial_xG\ast f=0,\\
v_t+(u+v+w)v_x+G\ast(vu_x+vw_x)+\partial_xG\ast g=0,\\
w_t+(u+v+w)w_x+G\ast(wu_x+wv_x)+\partial_xG\ast h=0,\\
 (u,v,w)|_{t=0}=(u_{0}(x),v_0(x),w_0(x)).\end{array}\right.
\end{equation}
We first recall the local well-posedness and blow-up phenomena which come from \cite{HLJ}.

\begin{lemma2} Assume the initial data $z_0=(u_0,v_0,w_0)\in H^{s}\times H^s\times H^s, s>\frac{3}{2}$.
Then there exists a unique strong solution $z=(u,v,w)$ to Eq.(2.3) and  a time $T=T(z_0)>0$, such that $$z(t,x)=(u,v,w)\in\mathcal {C}([0,T);H^{s}\times H^s\times H^s)\cap\mathcal {C}^1([0,T);
H^{s-1}\times H^{s-1}\times H^{s-1}).$$ Moreover, the solution $z(t,x)$ depends continuously on the
initial data $z_0$, i.e. the mapping $z_0 \rightarrow
z(\cdot,z_0)\; :$
$$H^{s}\times H^s\times H^s\rightarrow \mathcal {C}([0,T);H^{s}\times H^s\times H^s)\cap\mathcal
{C}^1([0,T); H^{s-1}\times H^{s-1}\times H^{s-1})$$ is continuous. Furthermore, the lifespan $T$ of solution $z(t,x)$ can be chosen independent of $s$.
\end{lemma2}
\begin{lemma2} Let $z_0=(u_0,v_0,w_0)\in H^{s}\times H^s\times H^s, s>\frac{3}{2}$, and T be the lifespan of solution $z=(u,v,w)$ to Eq.(2.3). Then it follows
for all $t\in [0,T)$ that
\begin{equation}\begin{split}E(t):=\|u(t)\|_{H^1}^2+\|v(t)\|_{H^1}^2&+\|w(t)\|_{H^1}^2\\&=\|u_0\|_{H^1}^2+\|v_0\|_{H^1}^2+\|w_0\|_{H^1}^2:=E_0,\end{split}\end{equation}
by the conservation law, we have
$$\|u(t)\|_{L^\infty}^2+\|v(t)\|_{L^\infty}^2+\|w(t)\|_{L^\infty}^2\leq\frac{1}{2}E_0.$$
Moreover, the solution $z=(u,v,w)$  blows up in finite time $T$ if and only if
\begin{equation}\liminf_{t\uparrow T}\inf_{x\in\mathbb{R}}\{u_x(x,t)\}=-\infty,\end{equation}
or\begin{equation}\liminf_{t\uparrow T}\inf_{x\in\mathbb{R}}\{v_x(x,t)\}=-\infty,\quad or\quad\liminf_{t\uparrow T}\inf_{x\in\mathbb{R}}\{w_x(x,t)\}=-\infty.\end{equation}
\end{lemma2}
\par
Next, we prove that there exists solutions to system (1.1) which do not exist globally in time. Comparing with the two  results of blow-up phenomena which are obtained  in \cite{HLJ}, we give another new wave-breaking solution.
\begin{theorem2}Let the initial data $z_0=(u_0,v_0,w_0)\in H^s\times H^s\times H^s,\;s>\frac{3}{2}$. Assume $T$ be the lifespan of solution $z=(u,v,w)$  to system (1.1). If the initial data $z_0$ satisfy
$$\int_{\mathbb{R}}(u_{0x}+v_{0x}+w_{0x})^3dx<-9E_0\sqrt{2E_0}.$$Then the corresponding solution $z(t,x)$ of system (1.1) blows up in finite time. Moreover, the lifespan $T$ is estimated above by
$$T\leq\frac{\sqrt{2E_0}}{3E_0} \log \left(\frac{Q(0)-9E_0\sqrt{2E_0}}{Q(0)+9E_0\sqrt{2E_0}}\right),$$
where $Q(0)=\int_{\mathbb{R}}(u_{0x}+v_{0x}+w_{0x})^3dx.$
\end{theorem2}
\textit{Proof.} Differentiating Eq.(2.2) with respective to $x$ variable, we have
\begin{equation}
\begin{split}
& u_{tx}=-(u+v+w)_xu_x-(u+v+w)u_{xx}-\partial_xG\ast(uv_x+uw_x)-\partial_x^2G\ast f,\\&
v_{tx}=-(u+v+w)_xv_x-(u+v+w)v_{xx}-\partial_xG\ast(vu_x+vw_x)-\partial_x^2G\ast g,\\&
w_{tx}-(u+v+w)_xw_x-(u+v+w)w_{xx}-\partial_xG\ast(wu_x+wv_x)-\partial_x^2G\ast h,
\end{split}
\end{equation}
where the functions $f,\;g$ and $h$ satisfy equality (2.2).\\
Thanks to (2.7), it follows that
\begin{equation}\begin{split}
\frac{\partial}{\partial t}(u_x&+v_x+w_x)^3=3(u_x+v_x+w_x)^2(u_{tx}+v_{tx}+w_{tx})\\&=-3(u_x+v_x+w_x)^4-3(u_x+v_x+w_x)^2[(u+v+w)(u+v+w)_{xx}\\&\quad+\partial_xG\ast(uv+uw+vw)_x
+\partial_x^2G\ast(f+g+h)].\end{split}
\end{equation}
Integrating Eq.(2.8) with respect to $x$ variable on $\mathbb{R}$ yields that
\begin{equation}\begin{split}
\frac{\partial}{\partial t}\int_{\mathbb{R}}&(u_x+v_x+w_x)^3dx=-2\int_{\mathbb{R}}(u_x+v_x+w_x)^4dx\\&-3\int_{\mathbb{R}}(u_x+v_x+w_x)^2[\partial_x^2G\ast(uv+uw+vw
+f+g+h)]dx.\end{split}
\end{equation}
By virtue of $\partial_x^2G\ast f=G\ast f-f$ we have
\begin{equation}\begin{split}-3\int_{\mathbb{R}}&(u_x+v_x+w_x)^2[\partial_x^2G\ast(uv+uw+vw
+f+g+h)]dx\\&=-3\int_{\mathbb{R}}(u_x+v_x+w_x)^2[G\ast(uv+uw+vw
+f+g+h)]dx\\&\quad+3\int_{\mathbb{R}}(u_x+v_x+w_x)^2[uv+uw+vw
+f+g+h]dx\\&=I+II.\end{split}\end{equation}
Note that
\begin{equation*}f+g+h=2(u^2+v^2+w^2)-\frac{1}{2}(u_x^2+v_x^2+w_x^2)+2(u_xv_x+u_xw_x+v_xw_x),
\end{equation*}
by Lemma 2.2 and $\|G\|_{L^\infty}\leq\frac{1}{2}$, the term $I$ can be bounded by
\begin{equation}\begin{split}I&\leq3\|(u_x+v_x+w_x)^2\|_{L^1}\|G\ast H\|_{L^\infty}\\&\leq\frac{9}{2}E_0\|H\|_{L^1}\\&
\leq\frac{27}{2}E_0^2,
\end{split}\end{equation}
where $$H=uv+uw+vw
+f+g+h.$$
We can deal with the term $II$ as follows
\begin{equation}\begin{split}II&=3\int_{\mathbb{R}}(u_x+v_x+w_x)^2[uv+uw+vw
+2(u^2+v^2+w^2)]dx\\&\quad+3\int_{\mathbb{R}}(u_x+v_x+w_x)^2[2(u_xv_x+u_xw_x+v_xw_x)-\frac{1}{2}(u_x^2+v_x^2+w_x^2)]dx\\&\leq9(\|u\|_{L^\infty}^2+
\|v\|_{L^\infty}^2+\|w\|_{L^\infty}^2)\int_{\mathbb{R}}(u_x+v_x+w_x)^2dx
\\&\quad+\frac{9}{2}\int_{\mathbb{R}}(u_x+v_x+w_x)^2(u_xv_x+u_xw_x+v_xw_x)]dx\\&
\leq\frac{27}{2}E_0^2+3\int_{\mathbb{R}}(u_x+v_x+w_x)^2[(u_xv_x+u_xw_x+v_xw_x)+\frac{1}{2}(u_x^2+v_x^2+w_x^2)]dx\\&
=\frac{27}{2}E_0^2+\frac{3}{2}\int_{\mathbb{R}}(u_x+v_x+w_x)^4dx.
\end{split}\end{equation}
Inserting (2.11) and (2.12) into (2.10). Then combining (2.9) with (2.10) to yield
\begin{equation}\frac{\partial}{\partial t}\int_{\mathbb{R}}(u_x+v_x+w_x)^3dx\leq-\frac{1}{2}\int_{\mathbb{R}}(u_x+v_x+w_x)^4dx+27E_0^2.
\end{equation}
By the following H\"{o}lder inequality
\begin{equation}\begin{split}
\left|\int_{\mathbb{R}}(u_x+v_x+w_x)^3dx\right|^2&\leq\int_{\mathbb{R}}(u_x+v_x+w_x)^4dx\int_{\mathbb{R}}(u_x+v_x+w_x)^2dx\\&\leq
3E_0\int_{\mathbb{R}}(u_x+v_x+w_x)^4dx,\end{split}\end{equation}
 and define$$Q(t)=\int_{\mathbb{R}}(u_x+v_x+w_x)^3dx.$$
 The inequality (2.13) is changed into
 \begin{equation}\begin{split}
\frac{\partial}{\partial t}Q(t)&\leq-\frac{1}{6E_0}Q^2(t)+27E_0^2\\&\leq-\frac{1}{6E_0}\left(Q(t)-9E_0\sqrt{2E_0}\right)\left(Q(t)+9E_0\sqrt{2E_0}\right).
 \end{split}
 \end{equation}
In view of the assumption $Q(0)<-9E_0\sqrt{2E_0}$ and (2.14), then  $\partial_tQ(t)<0$ and $Q(t)$ is a decreasing function, hence
 $$Q(t)<-9E_0\sqrt{2E_0}.$$
 By solving the inequality (2.15), one can easily check that
 \begin{equation}
 \frac{Q(0)+9E_0\sqrt{2E_0}}{Q(0)-9E_0\sqrt{2E_0}}e^{\frac{3}{2}\sqrt{2E_0}t}-1\leq\frac{18E_0\sqrt{2E_0}}{Q(t)-9E_0\sqrt{2E_0}}\leq0.
 \end{equation}
 Observing  that
 $$0<\frac{Q(0)+9E_0\sqrt{2E_0}}{Q(0)-9E_0\sqrt{2E_0}}<1.$$
 In view of (2.16), we can deduce that the lifespan $T$ of solution $z$ satisfies
 \begin{equation}
0<T\leq\frac{\sqrt{2}}{3\sqrt{E_0}} \log \left(\frac{Q(0)-9E_0\sqrt{2E_0}}{Q(0)+9E_0\sqrt{2E_0}}\right),
 \end{equation}
such that $\lim_{t\uparrow T}Q(t)=-\infty$. On the other hand
\begin{equation*}\begin{split}\left|\int_{\mathbb{R}}(u_x+v_x+w_x)^3dx\right|&\leq\|(u_x+v_x+w_x)\|_{L^\infty}\|(u_x+v_x+w_x)\|_{L^2}^2\\&
\leq3E_0\|(u_x+v_x+w_x)\|_{L^\infty}.
\end{split}
 \end{equation*}
which completes the proof Theorem 2.1. \hspace{\fill}$\square$

\section{The exponential decay of solution}
In this section, our aim is to establish the exponential decay of solution to system (1.1), before  stating precisely our main results, we first give two  important lemmas, which  will be continuously used in the  paper.
\begin{lemma3}(The  Gronwall Lemma) Let $f(t), g(t), h(t)$ be continuous functions on $\mathbb{R}^+$ such that
$$\partial_tf\leq g+hf.$$Then the following  inequality holds
\begin{equation}
f(t)\leq e^{\int_0^th(s)ds}f(0)+\int_0^tg(\tau)e^{\int_\tau^th(s)ds}d\tau.
\end{equation}
Moreover, if $f(t), g(t), h(t)$ are positive, and satisfy
$$f\leq g+\int_0^th(s)f(s)ds,$$then, we have
\begin{equation}f(t)\leq g(t)+\int_0^th(s)g(s)e^{\int_s^th(\tau)d\tau}ds.
\end{equation}
\end{lemma3}
\textit{Proof.} Applying $e^{-\int_0^th(s)ds}$ to  the inequality $\partial_tf\leq g+hf$ to immediately derive (3.1).
\par
Let
$$\mathcal {R}(t)=\int_0^th(s)f(s)ds.$$
Then the derivative $\mathcal {R}'$ satisfies
$$\mathcal {R}'(s)-h(s)\mathcal {R}(s)=h(s)(f(s)-\mathcal {R}(s))\leq h(s)g(s).$$
Consequently
$$\frac{d}{ds}\left\{\mathcal {R}(s)e^{\int_s^th(\tau)d\tau}\right\}\leq h(s)g(s)e^{\int_s^th(\tau)d\tau}.$$
Integrating on $[0,t]$ with respect to $s$ variable gives
$$\mathcal {R}(t)\leq\int_0^th(s)g(s)e^{\int_s^th(\tau)d\tau}ds.$$
Adding $g(t)$ on both sides of the above inequality to obtain (3.2).\hspace{\fill}$\square$

\begin{lemma3} Assume the function $G(x)=\frac{1}{2}e^{-|x|}$. Let the weighted function $J_N(x)$ be
\begin{equation}
J_N(x)=\left\{\begin{array}{ll}e^{\alpha N},\qquad &x\in (-\infty,-N),\\
e^{-\alpha x},\qquad &x\in[-N,0]\\
1,\qquad &x\in(0,\infty),
\end{array}\right.
\end{equation}
where $N\in \mathbb{Z}^+$.
If the constant $\alpha\in (0,1)$, then there exists some constant $C_0$, such that
\begin{equation}J_N(x)(G\ast ({J_N})^{-1})(x)\leq C_0\quad \text{and}\quad J_N(x)(\partial_xG\ast (J_N)^{-1})(x)\leq C_0.
\end{equation}
\end{lemma3}
\textit{Proof.} At first, we prove the first inequality. Note that
\begin{equation}\begin{split}
f(x):&=2J_N(x)(G\ast (J_N)^{-1})(x)=J_N(x)\int_{\mathbb{R}}\frac{e^{-|x-y|}}{J_N(y)}dy\\
&=J_N(x)\int_{x}^\infty\frac{e^{x-y}}{J_N(y)}dy+J_N(x)\int_{-\infty}^x\frac{e^{y-x}}{J_N(y)}dy.
\end{split}\end{equation}
Case 1: As $x<-N$, then we have
\begin{equation}\begin{split}
f(x)&=e^{\alpha N}\int_{(x,-N)\cup[-N,0]\cup(0,\infty)}\frac{e^{x-y}}{J_N(y)}dy+e^{\alpha N}\int_{-\infty}^x\frac{e^{y-x}}{J_N(y)}dy\\&
=e^{\alpha N}\left(\int_x^{-N}e^{x-y-\alpha N}dy+\int_{-N}^0e^{x-y+\alpha x}dy+\int_0^\infty e^{x-y}dy\right)\\&\qquad\qquad+e^{\alpha N}\int_{-\infty}^xe^{y-x-\alpha N}dy\\&=1-e^{x+N}+e^{(\alpha+1)(x+N)}-e^{\alpha(x+N)+x}+e^{\alpha N+x}\\&\leq3.
\end{split}\end{equation}
Case 2: If $x\in[-N,0]$,  we can derive
\begin{equation}\begin{split}
f(x)&=e^{-\alpha x}\int_{[x,0]\cup(0,\infty)}\frac{e^{x-y}}{J_N(y)}dy+e^{-\alpha x}\int_{(-\infty,-N)\cup(-N,x)}\frac{e^{y-x}}{J_N(y)}dy\\&
=e^{-\alpha x}\left(\int_x^{0}e^{x+(\alpha-1)y}dy+\int_0^\infty e^{x-y}dy\right)\\&\qquad+e^{-\alpha x}\left(\int_{-\infty}^{-N}e^{y-(x+\alpha N)}dy+\int_{-N}^xe^{(\alpha+1)y-x}dy\right)\\&=\frac{1}{1-\alpha}+e^{(1-\alpha)x}+e^{-(\alpha+1)(x+N)}+\frac{1}{\alpha+1}(1-e^{-(\alpha+1)(x+N)})\\
&\leq\frac{3-\alpha^2}{1-\alpha^2}.
\end{split}\end{equation}
Case 3: As $x\in (0,\infty)$, we can deal with it  as follows
\begin{equation}\begin{split}
f(x)&=\int_{x}^\infty e^{x-y}dy+\int_{(-\infty,-N)\cup[-N,0]\cup(0,x)}\frac{e^{y-x}}{J_N(y)}dy\\&
=1+e^{-(x+\alpha N+N)}+e^{(\alpha-1)x}-e^{(\alpha-1)x-N}+1-e^{-x}\\
&\leq3+e^{-(\alpha+1)N}-e^{-N}.
\end{split}\end{equation}
Let $2C_0\geq\frac{3-\alpha^2}{1-\alpha^2}$. Combining (3.5), (3.6), (3.7) with (3.8) to yield the first inequality. Similarly, the second inequality can be proved.\hspace{\fill}$\square$
\begin{remark3} If we define weighted function for $\alpha\in (0,1)$,
\begin{equation}
\varphi_N(x)=\left\{\begin{array}{ll}1,\qquad &x\in(-\infty,0),\\
e^{\alpha x},\qquad &x\in[0,N],\\
 e^{\alpha N},\qquad &x\in(N,\infty), \end{array}\right.
\end{equation}
where $N\in \mathbb{Z}^+$.  Then there exists some constant $C_0$, for all $N$, it follows  that
\begin{equation}\left\{\begin{array}{ll}
\varphi_N(x)(G\ast ({\varphi_N})^{-1})(x)\leq C_0,\\
\varphi_N(x)(\partial_xG\ast (\varphi_N)^{-1})(x)\leq C_0.
\end{array}\right.
\end{equation}
\end{remark3}
\par
 Next, as \cite{H-M, W-G}, we shall establish the exponential
decay of the strong solutions to Eq.(2.3), if
the initial data $z_0(x)$ decay at infinity.

\begin{theorem3} Let the initial data $z_0=(u_0,v_0,w_0)\in H^{s}\times H^{s}\times H^{s},s>\frac{3}{2}$ and
$T>0$. Suppose $z=(u,v,w)\in \mathcal {C}([0,T];H^s\times H^s\times H^s)$ is the
corresponding solution to Eq.(2.3) with the initial data $z_0$. If
there exists some $\alpha\in(0,1)$ such that
\begin{equation*}
\left\{\begin{array}{ll}|u_0(x)|,|v_0(x)|,|w_0(x)|\sim \mathcal{O}(e^{\alpha
x})\qquad \text{as}\quad x\downarrow-\infty,\\
\\ |u_{0x}(x)|,|v_{0x}(x)|, |w_{0x}(x)|\sim \mathcal{O}(e^{\alpha x})\qquad
\text{as}\quad x\downarrow-\infty, \end{array}\right.
\end{equation*} then, it follows that the solutions $z(t,x)$ satisfy
\begin{equation*}
\left\{\begin{array}{ll}|u(t,x)|,|v(t,x)|,|w(t,x)|\sim \mathcal{O}(e^{\alpha
x})\qquad \text{as}\quad x\downarrow -\infty,\\
\\ |u_{x}(t,x)|, |v_{x}(t,x)|, |w_{x}(t,x)|\sim \mathcal{O}(e^{\alpha x})\qquad
\text{as}\quad x\downarrow -\infty, \end{array}\right.
\end{equation*} uniformly in the interval
$[0,T]$.
\end{theorem3}
\textit{Proof.} For simplicity, let
$M=\sup_{t\in[0,T]}\{\|u(t)\|_{H^s}\}$,
using the Sobolev embedding theorem,
$\|u(t)\|_{L^{\infty}},\;\|u_x(t)\|_{L^{\infty}}\leq
M$.
\par
 Define a weighted function
\begin{equation}
J_N(x)=\left\{\begin{array}{ll}e^{\alpha N},\qquad &x\in (-\infty,-N),\\
e^{-\alpha x},\qquad &x\in[-N,0]\\
1,\qquad &x\in(0,\infty),
\end{array}\right.
\end{equation}

where $N\in \mathbb{Z}^+$. One can easily check that for all $N$
\begin{equation}
0\leq -J_N^{'}(x)\leq J_N(x) \quad a.e.\quad x\in
\mathbb{R}.
\end{equation}
Applying Eq.(2.3) by $J_N$ to deduce

\begin{equation}
(uJ_N)_t+(u+v+w)J_Nu_x+J_N[G\ast(uv_x+uw_x)]+J_N[\partial_xG\ast f]=0,
\end{equation}
\begin{equation}
(vJ_N)_t+(u+v+w)J_Nv_x+J_N[G\ast(vu_x+vw_x)]+J_N[\partial_xG\ast g]=0,
\end{equation}
and
\begin{equation}
(wJ_N)_t+(u+v+w)J_Nw_x+J_N[G\ast(wu_x+wv_x)]+J_N[\partial_xG\ast h]=0.
\end{equation}
Taking the scalar product of $(uJ_N)^{2p-1}$ and Eq.(3.13), integration by parts is given by the following equality
\begin{equation}\begin{split}
\int_{\mathbb{R}}(uJ_{N})_t(uJ_N)^{2p-1}&dx=-\int_{\mathbb{R}}J_N(u+v+w)u_x(uJ_N)^{2p-1}dx\\&-\int_{\mathbb{R}}\left(J_N[G\ast(uv_x+uw_x)]+J_N[\partial_xG\ast f]\right)(uJ_N)^{2p-1}dx.
\end{split}\end{equation}
Note that
$$\int_{\mathbb{R}}(uJ_N)_t(uJ_N)^{2p-1}dx=\|uJ_N\|_{L^{2p}}^{2p-1}\frac{d}{dt}
\|uJ_N\|_{L^{2p}},$$
\begin{equation*}\begin{split}
\int_{\mathbb{R}}J_N(u+v+w)u_x(uJ_N)^{2p-1}dx&\leq
\|u_x\|_{L^\infty}\|J_N(u+v+w)\|_{L^{2p}}\|uJ_N\|_{L^{2p}}^{2p-1}
\\&\leq
M\|(uJ_N,vJ_N,wJ_N)\|_{L^{2p}}\|uJ_N\|_{L^{2p}}^{2p-1},
\end{split}\end{equation*}
and
\begin{equation*}\begin{split}\int_{\mathbb{R}}(J_N&[G\ast(uv_x+uw_x)]+J_N[\partial_xG\ast f])(uJ_N)^{2p-1}dx\\&\leq\|uJ_N\|
_{L^{2p}}^{2p-1}(\|J_N[G\ast(uv_x+uw_x)]\|_{L^{2p}}+\|J_N[\partial_xG\ast f]\|_{L^{2p}}).\end{split}\end{equation*}
 In view of (3.16) and the above relations, we have
\begin{equation}
\frac{d}{dt} \|uJ_N\|_{L^{2p}}\leq
M\|(uJ_N,vJ_N,wJ_N)\|_{L^{2p}}+\|J_N[G\ast(uv_x+uw_x)]\|_{L^{2p}}+\|J_N[\partial_xG\ast f]\|_{L^{2p}},
\end{equation}
where $\|(uJ_N,vJ_N,wJ_N)\|_{L^{2p}}=\|uJ_N\|_{L^{2p}}+\|vJ_N\|_{L^{2p}}+\|wJ_N\|_{L^{2p}}$.
Since the function $\varphi\in L^{1}(\mathbb{R})\cap L^{\infty}(\mathbb{R})$ implies
$$\lim_{n\uparrow\infty}\|\varphi\|_{L^{n}}=\|f\|_{L^{\infty}}.$$
 Let $p\uparrow \infty$ in (3.17), it follows that
 \begin{equation}
\frac{d}{dt} \|uJ_N\|_{L^{\infty}}\leq
M\|(uJ_N,vJ_N,wJ_N)\|_{L^{\infty}}+\|J_N[G\ast(uv_x+uw_x)]\|_{L^{\infty}}+\|J_N[\partial_xG\ast f]\|_{L^{\infty}}.
\end{equation}
Similar to the process of (3.18), multiplying Eq.(3.14), Eq.(3.15) by $(vJ_N)^{2p-1}$ and $(wJ_N)^{2p-1}$, respectively, integrating the result on $\mathbb{R}$ with respect to $x$-variable, we end up with
\begin{equation}
\frac{d}{dt} \|vJ_N\|_{L^{\infty}}\leq
M\|(uJ_N,vJ_N,wJ_N)\|_{L^{\infty}}+\|J_N[G\ast(vu_x+vw_x)]\|_{L^{\infty}}+\|J_N[\partial_xG\ast g]\|_{L^{\infty}},
\end{equation}
and
\begin{equation}
\frac{d}{dt} \|wJ_N\|_{L^{\infty}}\leq
M\|(uJ_N,vJ_N,wJ_N)\|_{L^{\infty}}+\|J_N[G\ast(wu_x+wv_x)]\|_{L^{\infty}}+\|J_N[\partial_xG\ast h]\|_{L^{\infty}},
\end{equation}
where the functions $f,\;g$ and $h$ satisfy equality (2.2).

By virtue of Lemma 3.2, there exists a constant $C_0$ such that
\begin{equation}\begin{split}
|J_N[G\ast(uv_x+uw_x)]|&=\left|J_N(x)\int_{\mathbb{R}}\frac{e^{-|x-y|}}{2J_N(y)}J_N(y)(uv_x+uw_x)(y)
dy\right|\\&\leq C_0\|v_x+w_x\|_{L^\infty}\|uJ_N\|_{L^{\infty}}\\&\leq2C_0M\|
uJ_N\|_{L^{\infty}}.
\end{split}\end{equation}
and
\begin{equation}\begin{split}
|J_N[\partial_xG\ast f]|&=\left|J_N(x)\int_{\mathbb{R}}\frac{e^{-|x-y|}}{2J_N(y)}J_N(y)f(y)
dy\right|\\&\leq C\|(uJ_N,vJ_N,wJ_N,u_xJ_N,v_xJ_N,w_xJ_N)\|_{L^\infty}.
\end{split}\end{equation}
Similarly, we have
\begin{equation}\begin{split}
\|J_N&[G\ast(vu_x+vw_x)]\|_{L^{\infty}}+\|J_N[\partial_xG\ast g]\|_{L^{\infty}}\\&\leq C(\|
vJ_N\|_{L^{\infty}}+\|(uJ_N,vJ_N,wJ_N,u_xJ_N,v_xJ_N,w_xJ_N)\|_{L^\infty}),
\end{split}\end{equation}
\begin{equation}\begin{split}
\|J_N&[G\ast(wu_x+wv_x)]\|_{L^{\infty}}+\|J_N[\partial_xG\ast h]\|_{L^{\infty}}\\&\leq C(\|
wJ_N\|_{L^{\infty}}+\|(uJ_N,vJ_N,wJ_N,u_xJ_N,v_xJ_N,w_xJ_N)\|_{L^\infty}).
\end{split}\end{equation}

Add up (3.18), (3.19) with (3.20), plugging (3.21), (3.22), (3.23) and (3.24) into the inequality, by the Gronwall lemma to yield
\begin{equation}\begin{split}
\|(uJ_N,&vJ_N,wJ_N)\|_{L^{\infty}}\leq e^{Ct}\|(u_0J_N,v_0J_N,w_0J_N)\|_{L^{\infty}}\\&\quad+\int_0^t\|(uJ_N,vJ_N,wJ_N,u_xJ_N,v_xJ_N,w_xJ_N)\|_{L^\infty}(\tau)d\tau.
\end{split}\end{equation}
Differentiating  Eq.(2.3) with respect to $x$ variable, after
multiplying the result by $J_N$ it follows that
\begin{equation}
(u_xJ_N)_t+[(u+v+w)_xu_x+(u+v+w)u_{xx}]J_N+[\partial_xG\ast(uv_x+uw_x)+\partial_x^2G\ast f]J_N=0.
\end{equation}
\begin{equation}
(v_xJ_N)_t+[(u+v+w)_xv_x+(u+v+w)v_{xx}]J_N+[\partial_xG\ast(vu_x+vw_x)+\partial_x^2G\ast h]J_N=0.
\end{equation}
\begin{equation}
(w_xJ_N)_t+[(u+v+w)_xw_x+(u+v+w)w_{xx}]J_N+[\partial_xG\ast(wu_x+wv_x)+\partial_x^2G\ast h]J_N=0.
\end{equation}
Multiplying Eq.(3.26) by $(uJ_N)^{2p-1}$ with
$p\in\mathbb{Z}^+$ and integrating the result on $\mathbb{R}$ with
respect to $x$-variable, applying Holder's inequality, we have
\begin{equation}\begin{split}
\|u_xJ_N\|_{L^{2p}}^{2p-1}\frac{d}{dt}\|u_xJ_N\|_{L^{2p}}&\leq M\|(u_xJ_N,v_xJ_N,w_xJ_N)\|_{L^{2p}}\|u_xJ_N\|_{L^{2p}}^{2p-1}\\&\quad+\|[\partial_xG\ast(uv_x+uw_x)+\partial_x^2G\ast f]J_N\|_{L^{2p}}\|u_xJ_N\|_{L^{2p}}^{2p-1}\\&\quad-\int_{\mathbb{R}}J_N(u+v+w)u_{xx}(u_xJ_N)^{2p-1}dx.
\end{split}\end{equation}
Observing that
\begin{equation}\begin{split}
\int_{\mathbb{R}}&J_N(u+v+w)u_{xx}(u_xJ_N)^{2p-1}dx=\int_{\mathbb{R}}(u+v+w)(u_xJ_N)^{2p-1}[(u_{x}J_N)_x-u_xJ_N']dx\\&
=-\frac{1}{2p}\int_{\mathbb{R}}(u+v+w)_x(u_xJ_N)^{2p}dx-\int_{\mathbb{R}}(u+v+w)u_xJ_N'(u_xJ_N)^{2p-1}dx\\&
\leq C(\|u_xJ_N\|_{L^{2p}}^{2p}),
\end{split}\end{equation}
where we have applied $|J_N'|\leq J_N$.
Substituting (3.30) into (3.29), letting $p\uparrow \infty$ to obtain
\begin{equation}\begin{split}
\frac{d}{dt}\|u_xJ_N\|_{L^{\infty}}&\leq M\|(u_xJ_N,v_xJ_N,w_xJ_N)\|_{L^{\infty}}\\&\qquad+\|[\partial_xG\ast(uv_x+uw_x)+\partial_x^2G\ast f]J_N\|_{L^{\infty}}.
\end{split}\end{equation}
Multiplying Eq.(3.27) and Eq.(3.28) by $(v_xJ_N)^{2p-1}$ and $(w_xJ_N)^{2p-1}$, respectively, integrating the result on $\mathbb{R}$ with respect to $x$-variable, it follows that
\begin{equation}\begin{split}
\frac{d}{dt}\|v_xJ_N\|_{L^{\infty}}&\leq M\|(u_xJ_N,v_xJ_N,w_xJ_N)\|_{L^{\infty}}\\&\qquad+\|[\partial_xG\ast(vu_x+vw_x)+\partial_x^2G\ast g]J_N\|_{L^{\infty}},
\end{split}\end{equation}
and
\begin{equation}\begin{split}
\frac{d}{dt}\|w_xJ_N\|_{L^{\infty}}&\leq M\|(u_xJ_N,v_xJ_N,w_xJ_N)\|_{L^{\infty}}\\&\qquad+\|[\partial_xG\ast(wu_x+wv_x)+\partial_x^2G\ast h]J_N\|_{L^{\infty}}.
\end{split}\end{equation}
In view of Lemma 3.2, we can derive
\begin{equation}
\|[\partial_xG\ast(uv_x+uw_x)]J_N\|_{L^{\infty}}\leq C\|uJ_N\|_{L^\infty}.
\end{equation}
Thanks to $\partial_x^2G\ast f=G\ast f-f$, by Lemma 3.2 again to yield
\begin{equation}\begin{split}
\|[\partial_x^2G\ast f]J_N\|_{L^\infty}&\leq\|fJ_N\|_{L^\infty}+\|[G\ast f]J_N\|_{L^\infty}\\&
C\|(uJ_N,vJ_N,wJ_N,u_xJ_N,v_xJ_N,w_xJ_N)\|_{L^\infty}.
\end{split}\end{equation}
Consequently,
\begin{equation}
\|[\partial_xG\ast(vu_x+vw_x)+\partial_x^2G\ast g]J_N\|_{L^{\infty}}\leq C\|(uJ_N,vJ_N,wJ_N,u_xJ_N,v_xJ_N,w_xJ_N)\|_{L^\infty},
\end{equation}
\begin{equation}
\|[\partial_xG\ast(wu_x+wv_x)+\partial_x^2G\ast h]J_N\|_{L^{\infty}}\leq C\|(uJ_N,vJ_N,wJ_N,u_xJ_N,v_xJ_N,w_xJ_N)\|_{L^\infty}.
\end{equation}

Add up (3.31), (3.32) with (3.33), plugging (3.34), (3.35), (3.36) and (3.37) into the inequality. Then by virtue of Gronwall's inequality
implies
\begin{equation}\begin{split}
\|(u_xJ_N,&v_xJ_N,w_xJ_N)\|_{L^{\infty}}\leq e^{Ct}\|(u_{0x}J_N,v_{0x}J_N,w_{0x}J_N)\|_{L^{\infty}}\\&\quad+C\int_0^t\|(uJ_N,vJ_N,wJ_N,u_xJ_N,v_xJ_N,w_xJ_N)\|_{L^\infty}(\tau)d\tau.
\end{split}\end{equation}
where $C$ is constant depending only on $C_0,M$.
\par
 Let $$Z(t)=(\|(J_Nu(t),J_Nv(t),J_Nw(t))\|_{L^{\infty}}+\|(J_Nu_x(t),J_Nv_x(t),J_Nw_x(t))\|_{L^{\infty}}
).$$
 Applying Lemma 3.1 to (3.38), for all $t\in[0,T]$, there exists a constant
$\tilde{C}=\tilde{C}(C_0,M,T)$ such that
\begin{equation}\begin{split}
Z(t)\leq \tilde{C}Z(0)\leq\tilde{C}(\|(u_0,v_0,&w_0)\max(1,e^{-\alpha
x})\|_{L^{\infty}}\\&+\|(u_{0x},v_{0x},w_{0x}) \max(1,e^{-\alpha x})\|_{L^{\infty}}).
\end{split}\end{equation}
Letting $N\uparrow \infty$, from (3.19), for all $t\in[0,T]$, it follows for $x\leq0$ that
\begin{equation*}\begin{split}
(\|(u,v,w)&e^{-\alpha
x}\|_{L^{\infty}}+\|(u_x,v_x,w_x)e^{-\alpha
x}\|_{L^{\infty}})\\&
\leq\tilde{C}(\|(u_0,v_0,w_0)e^{-\alpha
x}\|_{L^{\infty}}+\|(u_{0x},v_{0x},w_{0x})e^{-\alpha x}\|_{L^{\infty}}),
\end{split}\end{equation*}
which obtains the desired result of Theorem 3.1. \hspace{\fill}$\square$
\par
If we choose the  weighted function $\varphi_N(x)$ for $\alpha\in (0,1)$ as
\begin{equation}
\varphi_N(x)=\left\{\begin{array}{ll}1,\qquad &x\in(-\infty,0),\\
e^{\alpha x},\qquad &x\in[0,N],\\
 e^{\alpha N},\qquad &x\in(N,\infty), \end{array}\right.
\end{equation}
where $N\in \mathbb{Z}^+$, then by virtue of Remark 3.1, by the method of proof of Theorem 3.1, we have the following result.
\begin{corollary3} Assume $z_0=(u_0,v_0,w_0)\in H^{s}\times H^{s}\times H^{s},s>\frac{3}{2}$ and
$T>0$. Suppose $z(t,x)=(u,v,w)\in \mathcal {C}([0,T];H^s\times H^s\times H^s)$ is the
corresponding solution to Eq.(2.3) with the initial data $z_0$. If
there exists some $\alpha\in(0,1)$ such that
\begin{equation*}
\left\{\begin{array}{ll}|u_0(x)|,|v_0(x)|,|w_0(x)|\sim \mathcal{O}(e^{-\alpha
x})\qquad \text{as}\quad x\uparrow\infty,\\
\\ |u_{0,x}(x)|,|v_{0,x}(x)|, |w_{0,x}(x)|\sim \mathcal{O}(e^{-\alpha x})\qquad
\text{as}\quad x\uparrow\infty, \end{array}\right.
\end{equation*} then the solutions $z$ satisfy
\begin{equation*}
\left\{\begin{array}{ll}|u(t,x)|,|v(t,x)|,|w(t,x)|\sim \mathcal{O}(e^{-\alpha
x})\qquad \text{as}\quad x\uparrow\infty,\\
\\ |u_{x}(t,x)|, |v_{x}(t,x)|, |w_{x}(t,x)|\sim \mathcal{O}(e^{-\alpha x})\qquad
\text{as}\quad x\uparrow \infty, \end{array}\right.
\end{equation*} uniformly in the interval
$[0,T]$.
\end{corollary3}
\begin{remark3}In fact, let $\alpha \in(0,1)$ and $j=0,1,2,\cdots,$ if the initial data $z_0$ satisfy $$(\partial_x^ju_0,\partial_x^jv_0,\partial_x^jw_0)\sim \mathcal{O}(e^{-\alpha |x|})\qquad
\text{as}\quad |x|\rightarrow \infty,$$ then the solutions $z$ to Eq.(2.3) satisfy
$$(\partial_x^ju,\partial_x^jv,\partial_x^jw)\sim \mathcal{O}(e^{-\alpha |x|})\qquad
\text{as}\quad |x|\rightarrow \infty.$$
\end{remark3}
Theorem 3.1 and Corollary 3.1 tell us that the solution $z$ can only decay as $e^{\alpha x}$ as $x\rightarrow-\infty$ and $e^{-\alpha x}$ as $x\rightarrow\infty$ for $\alpha\in(0,1)$. Whether the decay is optimal? the next result tell us some information.
\begin{theorem3}
 Given $z_0=(u_0,v_0,w_0)\in H^s\times H^s\times H^s,s\geq3$. Let $T=T(z_0)$ be the maximal existence time of the
solutions $z(t,x)=(u,v,w)$ to system (1.1) with the initial data $z_0$. If for
some $\lambda\geq0$ and $p\geq1$,
\begin{equation}
\|(m_0,n_0,l_0)e^{(1+\lambda)|x|}\|_{L^{2p}}\leq C,
\end{equation}
then we have for all $t\in[0,T)$ that
\begin{equation}
\|(m,n,l)e^{(1+\lambda)|x|}\|_{L^{2p}}\leq C.
\end{equation}
Moreover, if the initial data satisfy
\begin{equation}
\partial_x^ju_0,\partial_x^jv_0,\partial_x^jw_0\sim\mathcal
{O}(e^{-(1+\lambda)|x|})\quad\text{as}\quad|x|\rightarrow\infty,\;j=0,1,2,
\end{equation}
then for any $t\in[0,T)$, it follows that
$$(m,n,l)\sim\mathcal
{O}(e^{-(1+\lambda)|x|})\quad\text{as}\quad|x|\rightarrow\infty$$ and there exists some $\alpha\in(0,1)$ such that
$$\lim_{x\rightarrow+\infty}|(\partial_x^ju,\partial_x^jv,\partial_x^jw)e^{\alpha x}|\leq C,\;\lim_{x\rightarrow-\infty}
|(\partial_x^ju,\partial_x^jv,\partial_x^jw)e^{-\alpha x}|\leq C.$$
\end{theorem3}
\textit{Proof.} Multiplying system $(1.1)_1$ by $e^{(1+\lambda)|x|}$, after taking inner product with $(me^{(1+\lambda)|x|})^{2p-1}$ we
have
\begin{equation}\begin{split}
\|me^{(1+\lambda)|x|}\|_{L^{2p}}^{2p-1}&\frac{\partial}{\partial t}\|me^{(1+\lambda)|x|}\|_{L^{2p}}+\int_{\mathbb{R}}(m_xu+2m u_x)e^{(1+\lambda)|x|}(me^{(1+\lambda)|x|})^{2p-1}dx\\&+\int_{\mathbb{R}}e^{(1+\lambda)|x|}(mv+mw)_x(me^{(1+\lambda)|x|})^{2p-1}dx\\&\leq M\|(n,l)e^{(1+\lambda)|x|}\|_{L^{2p}}\|me^{(1+\lambda)|x|}\|_{L^{2p}}^{2p-1}.
\end{split}\end{equation}
Due to
\begin{equation}\begin{split}
\int_{\mathbb{R}}(m_xu+2m u_x)e^{(1+\lambda)|x|}&(me^{(1+\lambda)|x|})^{2p-1}dx=-\frac{2p-1}{2p}\int_{\mathbb{R}}u\partial_x(me^{(1+\lambda)|x|})^{2p}dx
\\&-\int_{\mathbb{R}}[-u_x+(1+\lambda)u\text{sgn}|x|](me^{(1+\lambda)|x|})^{2p}dx\\&\leq(\|u_x\|_{L^\infty}+(1+\lambda)\|u\|_{L^\infty})
\|me^{(1+\lambda)|x|}\|_{L^{2p}}^{2p}\\&\leq C\|me^{(1+\lambda)|x|}\|_{L^{2p}}^{2p},
\end{split}\end{equation}
and
\begin{equation}\begin{split}
\int_{\mathbb{R}}e^{(1+\lambda)|x|}(mv+mw)_x(me^{(1+\lambda)|x|})^{2p-1}dx
\leq C\|me^{(1+\lambda)|x|}\|_{L^{2p}}^{2p},
\end{split}\end{equation}
Combining (3.44), (3.45) with (3.46) to imply
\begin{equation}\begin{split}
\frac{\partial}{\partial t}\|me^{(1+\lambda)|x|}\|_{L^{2p}}\leq C\|(m,n,l)e^{(1+\lambda)|x|}\|_{L^{2p}}.
\end{split}\end{equation}
As the process of the estimation to (3.47), we deal with system $(1.1)_2$ and system $(1.1)_3$ is given by
\begin{equation}\begin{split}
\frac{\partial}{\partial t}\|ne^{(1+\lambda)|x|}\|_{L^{2p}}\leq C\|(m,n,l)e^{(1+\lambda)|x|}\|_{L^{2p}}.
\end{split}\end{equation}
\begin{equation}\begin{split}
\frac{\partial}{\partial t}\|le^{(1+\lambda)|x|}\|_{L^{2p}}\leq C\|(m,n,l)e^{(1+\lambda)|x|}\|_{L^{2p}}.
\end{split}\end{equation}
Add up (3.47), (3.48) with (3.49), then by the Gronwall inequality  yields that
\begin{equation}\begin{split}
\|(m,n,l)e^{(1+\lambda)|x|}\|_{L^{2p}}\leq C\|(m_0,n_0,l_0)e^{(1+\lambda)|x|}\|_{L^{2p}}.
\end{split}\end{equation}
By virtue of the assumption (3.41), it follows that (3.42).
In view of the assumption (3.43) to obtain
\begin{equation}
(m_0(x),n_0(x),l_0(x))\sim\mathcal
{O}(e^{-(1+\lambda)|x|})\quad\text{as}\quad|x|\uparrow\infty.
\end{equation}
Let $p\uparrow\infty$ in (3.50). Combining (3.50) with (3.51), let $|x|$ large enough, we have
$$(m,n,l)(t,x)\sim\mathcal
{O}(e^{-(1+\lambda)|x|})\quad\text{as}\quad|x|\uparrow\infty.$$
On the other hand, by virtue of (3.43), Theorem 3.1 and Corollary 3.1, we deduce for
any $\alpha\in(0,1)$ that
\begin{equation}
(\partial_x^ju,\partial_x^jv,\partial_x^jw)\sim\mathcal
{O}(e^{-\alpha|x|})\quad\text{as}\quad|x|\uparrow\infty\quad
j=0,\;1,2.
\end{equation}
This means that for all $t\in[0,T)$ and $j=0,\;1,\;2.$
$$\lim_{x\rightarrow+\infty}(\partial_x^ju,\partial_x^jv,\partial_x^jw)e^{\alpha x}\leq C,\;\lim_{x\rightarrow-\infty}
(\partial_x^ju,\partial_x^jv,\partial_x^jw)e^{-\alpha x}\leq C.$$
This completes the proof of Theorem 3.2. \hspace{\fill}$\square$

\begin{remark3} As long as the solution $z(t,x)$ exists, the result
of Theorem 3.2 tells us that the solutions $(z_,z_x)$ decay as $e^{\alpha x}$ as $x\rightarrow-\infty$ and $e^{-\alpha x}$ as $x\rightarrow\infty$ for $\alpha\in(0,1)$. However, the potential $(m,n,l)$ can decay as $e^{-(1+\lambda)|x|}$ as $|x|\rightarrow\infty$ for $\lambda\in(0,\infty)$.
\end{remark3}

\section{Traveling wave solutions}
In the subsection, we will establish a family of traveling
wave solutions to system (1.1).
\par
At first, we gives two important definitions and an useful lemma.
\begin{definition4} The solution $z(t,x)=(u,v,w)$ to system (1.1) is $x$-symmetric if there exists a function $b(t)\in
\mathcal {C}^{1}(\mathbb{R^{+}})$ such that $$z(t,x)=z(t,2b(t)-x),
\qquad\forall t\in[0,\infty),$$ for almost every $x\in \mathbb{R}$,
then the function $b(t)$ is  called the symmetric axis of $z(t,x)$.
\end{definition4}
\begin{definition4} Let $\mathcal {N}(\mathbb{R})=\{z:z=(u,v,w)\in
\mathcal {C}(\mathbb{R}^{+},H^{1}\times H^{1}\times H^1\}$. If for all $\phi\in \mathcal
{C}_{0}^{\infty}(\mathbb{R}^{+}\times\mathbb{R})$ and  $z(t,x)\in
\mathcal {N}(\mathbb{R})$ satisfy
\begin{equation}\begin{split}\langle u,&(1-\partial_{x}^{2})\phi_{t}\rangle+\left\langle
u_{x}(v+w),\phi_{xx}\right\rangle-\frac{1}{2}\left\langle u^2,
\phi_{xxx}\right\rangle+\\&\left\langle\frac{3}{2}u^2+\frac{1}{2}u_x^{2}+u(v+w)+u_x(v+w)_x+\frac{1}{2}(v^2+w^2-v_x^2-w_x^2),\phi_{x}\right\rangle=0,
\end{split}\end{equation}
\begin{equation}\begin{split}\langle v,&(1-\partial_{x}^{2})\phi_{t}\rangle+\left\langle
v_{x}(u+w),\phi_{xx}\right\rangle-\frac{1}{2}\left\langle v^2,
\phi_{xxx}\right\rangle+\\&\left\langle\frac{3}{2}v^2+\frac{1}{2}v_x^{2}+v(u+w)+v_x(u+w)_x+\frac{1}{2}(u^2+w^2-u_x^2-w_x^2),\phi_{x}\right\rangle=0,
\end{split}\end{equation}
\begin{equation}\begin{split}\langle w,&(1-\partial_{x}^{2})\phi_{t}\rangle+\left\langle
w_{x}(u+v),\phi_{xx}\right\rangle-\frac{1}{2}\left\langle w^2,
\phi_{xxx}\right\rangle+\\&\left\langle\frac{3}{2}w^2+\frac{1}{2}w_x^{2}+w(u+v)+w_x(u+v)_x+\frac{1}{2}(u^2+v^2-u_x^2-v_x^2),\phi_{x}\right\rangle=0,
\end{split}\end{equation}
then $z(t,x)$ is
a weak solution to system (1.1), where  $\langle\cdot,\cdot\rangle$ denotes the distributions on (t,x).
\end{definition4}
\begin{lemma4} Assume that $Z(x)=(U,V,W)\in \mathcal {N}(\mathbb{R})$ and satisfies
\begin{equation}\begin{split}&\langle -cU,(1-\partial_{x}^{2})\phi_{x}\rangle+\left\langle
U_{x}(V+W),\phi_{xx}\right\rangle-\frac{1}{2}\left\langle U^2,
\phi_{xxx}\right\rangle+\\&\left\langle\frac{3}{2}U^2+\frac{1}{2}U_x^{2}+U(V+W)+U_x(V+W)_x+\frac{1}{2}(V^2+W^2-V_x^2-W_x^2),\phi_{x}\right\rangle=0,
\end{split}\end{equation}
\begin{equation}\begin{split}&\langle -cV,(1-\partial_{x}^{2})\phi_{x}\rangle+\left\langle
V_{x}(U+W),\phi_{xx}\right\rangle-\frac{1}{2}\left\langle V^2,
\phi_{xxx}\right\rangle+\\&\left\langle\frac{3}{2}V^2+\frac{1}{2}V_x^{2}+V(U+W)+V_x(U+W)_x+\frac{1}{2}(U^2+W^2-U_x^2-W_x^2),\phi_{x}\right\rangle=0,
\end{split}\end{equation}
\begin{equation}\begin{split}&\langle -cW,(1-\partial_{x}^{2})\phi_{x}\rangle+\left\langle
W_{x}(U+V),\phi_{xx}\right\rangle-\frac{1}{2}\left\langle W^2,
\phi_{xxx}\right\rangle+\\&\left\langle\frac{3}{2}W^2+\frac{1}{2}W_x^{2}+W(U+V)+W_x(U+V)_x+\frac{1}{2}(U^2+V^2-U_x^2-V_x^2),\phi_{x}\right\rangle=0,
\end{split}\end{equation}
for all $\phi\in \mathcal {C}_{0}^{\infty}(\mathbb{R}^+\times\mathbb{R})$. Then the
function $z$  is given by
\begin{equation} z(t,x)=Z(x-c(t-t_{0}))
\end{equation}
is a weak solution of system (1.1), for any fixed $t_{0}\in\mathbb{R}^+$.
\end{lemma4}
\textit{Proof.}  Since $\mathcal {C}_{0}^{\infty}(\mathbb{R}^{+}\times\mathbb{R})$ is dense in
$\mathcal {C}_{0}^{1}(\mathbb{R}^{+},\mathcal
{C}_{0}^{3}(\mathbb{R}))$, by the density argument, we only consider
the test functions belongs to $\mathcal
{C}_{0}^{1}(\mathbb{R}^{+},\mathcal {C}_{0}^{3}(\mathbb{R}))$. Without loss of generality, let $t_{0}=0$.
Choosing  $\psi\in\mathcal
{C}_{0}^{1}(\mathbb{R}^{+},\mathcal {C}_{0}^{3}(\mathbb{R}))$, let
 $\psi_{c}=\psi(t,x+ct)$, we can derive
\begin{equation}
\partial_x(\psi_{c})=(\psi_{x})_{c},\qquad \text{and}\quad
\partial_t(\psi_{c})=(\psi_{t})_{c}+c(\psi_{x})_{c}.
\end{equation}
Assume $z(t,x)=Z(x-c(t-t_{0}))$. It is easy to check
\begin{equation}\left\{\begin{array}{ll}\langle
u,\psi\rangle=\langle U,\psi_{c}\rangle,\;\langle
u^{2},\psi\rangle=\langle U^{2},\psi_{c}\rangle,\\
\\\langle u_{x}v,\psi\rangle=\langle
U_{x}V,\psi_{c}\rangle,\;\langle u_{x}^{2},\psi\rangle=\langle
U_{x}^{2},\psi_{c}\rangle,\end{array}\right.
\end{equation}
where $Z=Z(x)$. In view of (4.8) and (4.9), it follows that
\begin{equation}\begin{split}
&\left\langle u,(1-\partial_{x}^{2})\psi_{t}\right\rangle
+\left\langle u_{x}(v+w),\psi_{xx}\right\rangle-\frac{1}{2}\left\langle u^2,
\psi_{xxx}\right\rangle\\&= \left\langle
U,(1-\partial_{x}^{2})(\psi_{t})_{c} \right\rangle+\left\langle U_{x}(V+W),(\psi_{xx})_c\right\rangle-\frac{1}{2}\left\langle U^2,
(\psi_{xxx})_c\right\rangle\\&
=\left\langle
U,(1-\partial_{x}^{2})(\partial_{t}\psi_{c}-c\partial_{x}\psi_{c})\right\rangle+
\left\langle U_{x}(V+W),\partial_x^2\psi_c\right\rangle-\frac{1}{2}\left\langle U^2,
\partial_x^3\psi_c\right\rangle,
\end{split}\end{equation}
and
\begin{equation}\begin{split}
\frac{1}{2}\langle 3u^2+u_x^{2}&+2u(v+w)+2u_x(v+w)_x+(v^2+w^2-v_x^2-w_x^2),\psi_{x}\rangle
\\&=\frac{1}{2}
\left\langle 3U^2+U_x^{2}+(V^2+W^2-V_x^2-W_x^2),\partial_x\psi_{c}\right\rangle\\&\qquad
+\left\langle U(V+W)+U_x(V+W)_x,\partial_x\psi_{c}\right\rangle.
\end{split}\end{equation}
Note that $Z$  only depends on $x$ variable, let $T$ large enough such that
it does not belong to the support of $\psi_{c}$, consequently
\begin{equation}\begin{split}
\left\langle U,(1-\partial_{x}^{2})\partial_{t}\psi_{c}\right
\rangle&=\int_{\mathbb{R}}U(x)\int_{\mathbb{R}^{+}}\partial_{t}(1-\partial_{x}^{2})\psi_{c}dtdx\\
&=\int_{\mathbb{R}}U(x)[(1-\partial_{x}^{2})\psi_{c}(T,x)-(1-\partial_{x}^{2})\psi_{c}(0,x)]dx\\&=0.
\end{split}\end{equation}
Combining (4.10), (4.11) with (4.12), it follows that
\begin{equation*}\begin{split}
&\left\langle u,(1-\partial_{x}^{2})\psi_{t}\right\rangle
+\left\langle u_{x}(v+w),\psi_{xx}\right\rangle-\frac{1}{2}\left\langle u^2,
\psi_{xxx}\right\rangle+\\&\quad\frac{1}{2}\langle 3u^2+u_x^{2}+2u(v+w)+2u_x(v+w)_x+(v^2+w^2-v_x^2-w_x^2),\psi_{x}\rangle\\&=
\langle -cU,(1-\partial_{x}^{2})\partial_x\psi_{c}\rangle+\left\langle
U_{x}(V+W),\partial_x^2\psi_{c}\right\rangle-\frac{1}{2}\left\langle U^2,
\partial_x^3\psi_{c}\right\rangle+\\&\quad\left\langle\frac{3}{2}U^2+\frac{1}{2}U_x^{2}+U(V+W)+U_x(V+W)_x+\frac{1}{2}(V^2+W^2-V_x^2-W_x^2),\partial_x\psi_{c}\right\rangle
=0,
\end{split}\end{equation*}
where we have applied (4.4) with $\phi(x)=\psi_{c}(t,x)$.
Therefore $u(t,x)=U(x-c(t-t_{0}))$ is a weak solution of system $(1.1)_1$.
Similarly, thanks to (4.5) and (4.6), we imply that $v(t,x)=V(x-c(t-t_{0}))$, $w(t,x)=W(x-c(t-t_{0}))$ is weak solutions to system $(1.1)_2$, $(1.1)_3$ respectively. This completes the proof of Lemma 4.1. \hspace{\fill}$\square$
\par
Finally, we state the main result in this subsection.
\begin{theorem4}
Assume $z(t,x)$ be $x$-symmetric. If  $z=(u,v,w)$ is a unique weak
solution of system (1.1), then $z(t,x)$ is a traveling wave.
\end{theorem4}
\textit{Proof.} It is necessary to consider the test function $\varphi\in
\mathcal {C}_{0}^{1}(\mathbb{R}^{+},\mathcal{C}_{0}^{3}(\mathbb{R}))$.
Let$$\varphi_{b}(t,x)=\varphi(t,2b(t)-x),\qquad b(t)\in \mathcal
{C}^{1}(\mathbb{R}).$$ One can easily check that
$(\varphi_{b})_{b}=\varphi$ and
\begin{equation}
\left\{\begin{array}{ll}
\partial_{t}\varphi_{b}=(\partial_{t}\varphi)_{b}+2\dot{b}(\partial_{x}\varphi)_{b},\\
\\
\partial_{x}u_{b}=-(\partial_{x}u)_{b},\;\partial_{x}\varphi_{b}
=-(\partial_{x}\varphi)_{b}.
\end{array}\right.
\end{equation}
Moreover,
\begin{equation}\left\{\begin{array}{ll}\langle
u_{b},\varphi\rangle=\langle u,\varphi_{b}\rangle,\;\langle
u_{b}^{2},\varphi\rangle=\langle u^{2},\varphi_{b}\rangle,\\
\\\langle v_{b}\partial_{x}u_{b},\varphi\rangle=-\langle
v\partial_{x}u,\varphi_{b}\rangle,\;\langle
(\partial_{x}u_{b})^{2},\varphi\rangle=\langle
(\partial_{x}u)^{2},\varphi_{b}\rangle,\end{array}\right.
\end{equation}
where $\dot{b}$ denotes the time derivative of $b$. \par Since $z$
is $x$-symmetric, in view  of (4.13) and (4.14),  we imply that
\begin{equation}\begin{split}
&\left\langle u,(1-\partial_{x}^{2})\varphi_{t}\right\rangle
+\left\langle u_{x}(v+w),\varphi_{xx}\right\rangle-\frac{1}{2}\left\langle u^2,
\varphi_{xxx}\right\rangle\\&= \left\langle
u,((1-\partial_{x}^{2})\partial_t\varphi)_{b}\right\rangle-\left\langle u_{x}(v+w),(\partial_x^2\varphi)_b\right\rangle-\frac{1}{2}\left\langle u^2,
(\partial_x^3\varphi)_b\right\rangle\\&
=\left\langle
u,(1-\partial_{x}^{2})(\partial_{t}\varphi_{b}+2\dot{b}\partial_{x}\varphi_{b})\right\rangle-
\left\langle u_{x}(v+w),\partial_x^2\varphi_b\right\rangle+\frac{1}{2}\left\langle u^2,
\partial_x^3\varphi_b\right\rangle,
\end{split}\end{equation}
and
\begin{equation}\begin{split}
\frac{1}{2}&\langle 3u^2+u_x^{2}+2u(v+w)+2u_x(v+w)_x+(v^2+w^2-v_x^2-w_x^2),\varphi_{x}\rangle
\\&=-\frac{1}{2}\langle 3u^2+u_x^{2}+2u(v+w)+2u_x(v+w)_x+(v^2+w^2-v_x^2-w_x^2),\partial_x\varphi_{b}\rangle.
\end{split}\end{equation}
Add up (4.15) with (4.16), by (4.1) we have
\begin{equation}\begin{split}&\left\langle u,(1-\partial_{x}^{2})\varphi_{t}\right\rangle
+\left\langle u_{x}(v+w),\varphi_{xx}\right\rangle-\frac{1}{2}\left\langle u^2,
\varphi_{xxx}\right\rangle+\\&\quad\frac{1}{2}\langle 3u^2+u_x^{2}+2u(v+w)+2u_x(v+w)_x+(v^2+w^2-v_x^2-w_x^2),\varphi_{x}\rangle\\&=\left\langle
u,(1-\partial_{x}^{2})(\partial_{t}\varphi_{b}+2\dot{b}\partial_{x}\varphi_{b})\right\rangle-
\left\langle u_{x}(v+w),\partial_x^2\varphi_b\right\rangle+\frac{1}{2}\left\langle u^2,
\partial_x^3\varphi_b\right\rangle-\\&\quad\frac{1}{2}\langle 3u^2+u_x^{2}+2u(v+w)+2u_x(v+w)_x+(v^2+w^2-v_x^2-w_x^2),\partial_x\varphi_{b}\rangle\\&
=0.
\end{split}\end{equation}
Thus, taking place $\varphi$ by $\varphi_{b}$ in (4.17), due to
$(\varphi_{b})_{b}=\varphi$, it follows that
\begin{equation}\begin{split}-\frac{1}{2}\langle &3u^2+u_x^{2}+2u(v+w)+2u_x(v+w)_x+(v^2+w^2-v_x^2-w_x^2),\partial_x\varphi\rangle\\&\quad+\left\langle
u,(1-\partial_{x}^{2})(\partial_{t}\varphi+2\dot{b}\partial_{x}\varphi)\right\rangle-
\left\langle u_{x}(v+w),\partial_x^2\varphi\right\rangle+\frac{1}{2}\left\langle u^2,
\partial_x^3\varphi\right\rangle\\&=0.
\end{split}\end{equation}
Subtracting (4.18) from (4.17) to derive
\begin{equation}\begin{split}
&\langle u,2\dot{b}(1-\partial_{x}^{2})\partial_{x}\varphi\rangle-2
\left\langle u_{x}(v+w),\partial_x^2\varphi\right\rangle+\left\langle u^2,
\partial_x^3\varphi\right\rangle-\\&\langle 3u^2+u_x^{2}+2u(v+w)+2u_x(v+w)_x+(v^2+w^2-v_x^2-w_x^2),\partial_x\varphi\rangle=0.
\end{split}\end{equation}
 If we choose
$\varphi_{\varepsilon}(t,x)=\psi(x)\varrho_{\varepsilon}(t)$ in (4.19), where $\psi\in \mathcal {C}_{0}^{\infty}(\mathbb{R})$ and
$\varrho_{\varepsilon}\in \mathcal {C}_{0}^{\infty}(\mathbb{R}^{+})$ is
a mollifier with the property that $\varrho_{\varepsilon}\rightarrow
\delta(t-t_{0})$, the Dirac mass at $t_{0}$, as
$\varepsilon\rightarrow0$. This implies that
\begin{equation}\begin{split}
\int_{\mathbb{R}}&(1-\partial_{x}^{2})\partial_{x}\psi\int_{\mathbb{R}^{+}}\dot{b}u
\varrho_{\varepsilon}(t)dtdx+\frac{1}{2}\int_{\mathbb{R}}\partial_{x}^3\psi\int_{\mathbb{R}^{+}}u^{2}
\varrho_{\varepsilon}(t)dtdx\\&-\int_{\mathbb{R}}\partial_x^2\psi\int_{\mathbb{R}^{+}}u_x(v+w)\varrho_{\varepsilon}(t)dtdx-\frac{1}{2}\int_{\mathbb{R}}
\partial_{x}\psi\int_{\mathbb{R}
^{+}}[2u_x(v+w)_x]\varrho_{\varepsilon}(t)dtdx\\
&-\frac{1}{2}\int_{\mathbb{R}}
\partial_{x}\psi\int_{\mathbb{R}
^{+}}[3u^2+u_x^{2}+2u(v+w)+(v^2+w^2-v_x^2-w_x^2)]\varrho_{\varepsilon}(t)dtdx=0.
\end{split}\end{equation}
Note that
\begin{equation*}
\lim_{\varepsilon\rightarrow0}\int_{\mathbb{R}_{+}}\dot{b}u
\rho_{\varepsilon}(t)dt=\dot{b}(t_{0})u(t_{0},x)\quad \text{in}\;
L^{2}(\mathbb{R}),
\end{equation*}
 Therefore, letting
$\varepsilon\rightarrow0$, (4.20) is given by
\begin{equation}\begin{split}
&\langle-\dot{b}(t_0)u(t_0,x),(1-\partial_{x}^{2})\partial_{x}\psi\rangle
+\langle u_x(v+w)(t_{0},x),\partial_x^2\psi\rangle-\frac{1}{2}\langle u^2(t_{0},x),\partial_x^3\psi\rangle
+\\&\frac{1}{2}\langle[3u^2+u_x^{2}+2u(v+w)+2u_x(v+w)_x+(v^2+w^2-v_x^2-w_x^2)](t_{0},x),\partial_x\psi\rangle\\&
\qquad=0.
\end{split}\end{equation}
Hence set $c=\dot{b}(t_{0})$, we prove that $u(t_{0},x)$ satisfies (4.4). As the process of (4.21), $v(t_{0},x)$, $w(t_{0},x)$ is the solution to (4.5), (4.6) respectively.
 By virtue of Lemma 4.1, $\tilde{z}(t,x)=z(t_{0},x-\dot{b}(t_0)(t-t_0))$ is a traveling wave
solution of system (1.1). In view of $\tilde{z}(t_0,x)=z(t_{0},x)$ and the
uniqueness of the solution of system(1.1), it follows that
$\tilde{z}(t,x)=z(t,x)$, for any $t>0$,
which conclude the proof of Theorem 4.1. \hspace{\fill}$\square$


\section*{Acknowledgments} This work was partially supported by CPSF (Grant No.: 2013T60086) and  NSFC (Grant No.: 11401122).
The author thanks the professor Boling Guo for his helpful discussions and
constructive suggestions and would like to say thanks to Pro. Qiaoyi Hu for sending several of her
papers to the author.


\end{document}